\documentclass[12pt]{article}
\usepackage{bm,latexsym,amsmath,amsthm,amssymb}
\usepackage{graphicx,indentfirst}
\usepackage{psfrag}
\usepackage{graphicx,amssymb}
\usepackage{latexsym,bm}
\begin{document}
\title{\bf  Global attractivity of  almost periodic solutions for  competitive Lotka-Volterra
diffusion system 
 \thanks{This work was supported by the
National Natural Science Foundation of China (Grant Nos. 11271312,
 11261058), the China Postdoctoral Science Foundation
(Grant Nos. 20110491750), the Natural Science Foundation of Xinjiang
(Grant Nos. 2011211B08).}}
\author{ Ahmadjan Muhammadhaji, Zhidong Teng\thanks{Corresponding author.
\  Tel/Fax: \ +86 991 8585505,
 Email:  zhidong@xju.edu.cn(Z. Teng), ahmatjanam@yahoo.com.cn(A. Muhammadhaji).}, Mehbuba Rehim \\
\small College of Mathematics and Systems Science,
Xinjiang University\\
\small Xinjiang, Urumqi 830046, P.R.China}
\date{\small{Received: 23-12-2012} / \small{Accepted: 23-01-2013/} \small{\textbf{Accepted for publication in AMV}}
}
\maketitle \baselineskip 18pt
\begin{center}
\begin{minipage}{130mm}
\textbf{Abstract}:\quad In this paper, two competitive
Lotka-Volterra populations in the two-patch-system with diffusion
are considered. Each of the two spiecies can diffuse indepently and
discretely between its in intrapatch and interpatch. By means of
constructing Liapunov function, under moderate condition, the system
has a unique almost periodic solution and
 which is  asymptotically stable and globally attractive  .

\end{minipage}
\end{center}

\section*{1. Introduction}
Diffusion is a ubiquitous phenomenon in the real world. It is
population pressure due to the mutual interference between the
individuals, describing the migration of species to avoid crowds. It
is important for us to understand the dynamics of populations of
nature, and the basic and important studied questions for the
dynamics of populations are the persistence, permanence and
extinction of species, global stability of systems and the existence
of positive periodic solutions,
 positive almost periodic solutions and asymptotically periodic solutions, etc.
Recently, many scholars have paid attention to the non-autonomous
Lotka-Volterra population models with diffusion.
There exists an extensive literature concerning the study of global stability and the existence of positive
 periodic solutions,
positive almost periodic solutions and asymptotically periodic solutions of Lotka-Volterra system
with diffusion and periodic parameters, see the monographs [1-3,5-19] and the references cited therein.

In [12], the authors studied the following nonautonomous Lotka-Volterra almost periodic cooperative
systems with diffusion
$$\left\{
    \begin{array}{ll}
      \dot{x}_{1}={x_{1}[r_{1}(t)-a_{11}(t)x_{1}+a_{12}(t)y_{1}]+D_{1}(t)(x_{2}-x_{1}),}\vspace{1.5mm}  \\
      \dot{y}_{1}={y_{1}[r_{2}(t)+a_{21}(t)x_{1}-a_{22}(t)y_{1}]+D_{2}(t)(y_{2}-y_{1}),}\vspace{1.5mm}  \\
      \dot{x}_{2}={x_{2}[s_{1}(t)-b_{11}(t)x_{2}+b_{12}(t)y_{2}]+D_{1}(t)(x_{1}-x_{2}),}\vspace{1.5mm}  \\
      \dot{y}_{2}={y_{2}[s_{2}(t)+b_{21}(t)x_{2}-b_{22}(t)y_{2}]+D_{2}(t)(y_{1}-y_{2}),}
\end{array}
  \right.\eqno(1.1)
$$
By means of constructing Liapunov function and under appropriate conditions,
the sufficient conditions on the existence of
 a unique almost periodic solution and its global asymptotic stability are
 established for  system (1.1). Based on the system (1.1), in [13], the authors
 generalized almost periodic system (1.1) into asymptotically periodic systems, under suitable conditions,
 the authors obtained that
asymptotically periodic systems have a unique solution which is globally asymptotically stable.

In [11], the authors studied the following nonautonomous
Lotka-Volterra  periodic competitive systems with diffusion
$$\left\{
    \begin{array}{ll}
      \dot{x}_{1}={x_{1}[r_{1}(t)-a_{11}(t)x_{1}-a_{12}(t)y_{1}]+D_{1}(t)(x_{2}-x_{1}),}\vspace{1.5mm}  \\
      \dot{y}_{1}={y_{1}[r_{2}(t)-a_{21}(t)x_{1}-a_{22}(t)y_{1}]+D_{2}(t)(y_{2}-y_{1}),} \vspace{1.5mm} \\
      \dot{x}_{2}={x_{2}[s_{1}(t)-b_{11}(t)x_{2}-b_{12}(t)y_{2}]+D_{1}(t)(x_{1}-x_{2}),}\vspace{1.5mm}  \\
      \dot{y}_{2}={y_{2}[s_{2}(t)-b_{21}(t)x_{2}-b_{22}(t)y_{2}]+D_{2}(t)(y_{1}-y_{2}),}
\end{array}
  \right.\eqno(1.2)
$$
By using of Brouwer fixed point theorem and
constructing a suitable Liapunov function, under some appropriate conditions,
 the authors obtained that the system
has a unique periodic solution which is globally stable.

Motivated by the works [12] and [13] of Wei and Wang, by using of
Liapunov method used in [12,13,14], we generalize system (1.2) into
almost periodic system, under suitable conditions, we obtain that
the system has a unique almost periodic solution which is
asymptotically stable and globally attractive.

The organization of this paper is as follows. In the next section we
will present some basic assumptions, notations and Lemmas . In
section 3, conditions for the almost periodic solution and
asymptotic stability are considered. In section 4, conditions for
the global attractivity are considered. In the final section, one
example is given to illustrate that our main results are applicable.

\section*{ 2.  Preliminaries}
In system (1.2), we have that
 $x_{1}(t), y_{1}(t)$ are the density of two competitive
species at time t at the first patch, $x_{2}(t) , y_{2}(t)$ are the
density of two competitive species at time t at the second patch,
$r_{i}(t)$ and $s_{i}(t)$ are intrinsic growth rate of two
competitive species at the first and second patch respectively,
$a_{ii}(t)$ and $b_{ii}(t)$ are intrapatch restriction density of
each species in two-patch-system, $a_{ij}(t), b_{ij}(t)(i\neq j)$
are competitive coefficients between two species, $D_{i}(t)$ are
diffusion coefficients. In this paper, we always assume that system
(1.2) satisfies the following assumption

$(H_1)$  $r_{i}(t),s_{i}(t),a_{ij}(t),b_{ij}(t)$ and $D_{i}(t)$ are
nonnegative continuous bounded almost periodic functions $(i, j = 1,
2)$.

From the viewpoint of mathematical biology, in this paper for system
(1.2) we only consider the solution with the following initial
condition
$$
x_i(t)=\phi_i(t),\ y_i(t)=\varphi_i(t) \quad\mbox{for all}\quad
t\in[0,+\infty),\; i=1,2
$$
where $\phi_i(t), \varphi_i(t)\;(i=1,2)$ are nonnegative continuous
functions defined on $[0,+\infty)$ satisfying $\phi_i(0)>0,
\varphi_i(0)>0 \;(i=1,2)$.

For a continuous and bounded function $f(t)$ defined on
$[0,+\infty)$, we define $f^{L}=\inf_{t\in[0,+\infty)}\{f(t)\}$ and
$ f^{M}=\sup_{t\in[0,+\infty)}\{f(t)\} $.

Now, we present some useful lemmas.

\textbf{Lemma 2.1. [11]} $R^4_{+0}=\{(x_1,y_1,x_2,y_2)\in
R^4|x_i\geq0,\ y_i\geq0,\ \ (i=1,2) \}$ is the positive invariance
set with respect to the system (1.2).\\

\textbf{Lemma 2.2. [11]} If the following inequalities hold
$$
D^M_1 < r^L_1,\quad  D^M_2 < r^L_2,\quad  D^M_1<s^L_1, \quad D^M_2 <
s^L_2,
$$
then there exists a compact region  which has a positive distance
from the coordinate hyperplane and it attracts all the solutions of
the system (1.2) with positive initial values.\\

\textbf{Lemma 2.3. [12]} Let $D$ be an open set of $R^4_+$, function
$V(t,x,y)$ be defined on the region $R_+\times D\times D$ or
  $R_+\times R_+^4\times R_+^4$, it satisfies that:

(i)$ a(||x - y||)\leq V (t, x, y)\leq b(||x - y||),$ where $a(r)$
and $b(r)$ are continuous increasing positive functions;

 (ii)$ ||V
(t,x_{1}, y_{1})- V (t,x_{2}, y_{2})||\leq k(||x_{1} - x_{2}|| +
||y_{1} - y_{2}||)$, $k>0$ is a constant;

(iii)$V'(t, x, y)\leq cV(t, x, y)$, where $c > 0$ is a constant.\\
\noindent Further, let the solution of system (1.2) lie in compact
set $\Omega$ for all $t\geq t_{0}>0,\ \Omega\in D$, then system
(1.2) has a unique almost periodic solution $z(t)$ in $D$, $z(t)$ lies
in $\Omega$, and it is uniformly asymptotically stable.\\

\textbf{Lemma 2.4. [4]} Let $f$ be a nonnegative function defined
on$[0,\infty)$ such that $f$ is integrable on $[0,\infty)$ and
uniformly continuous on$ [0,\infty)$.
Then $ \lim_{t\rightarrow\infty}f(t)=0$ .\\

F.Wei et al. obtained in [11] that the system (1.2) had a bounded
closed and convex set
$$
\Omega=\bigg\{z:z\in R_+^4, S(z)\leq\beta,x_i^L\leq x_1\leq x_i^M,y_i^L\leq y_1\leq y_i^M(i=1,2),
R_+=[0,+\infty)\bigg\}.
$$
where $S(z),\beta, x_i^L, x_i^M, y_i^L, y_i^M$ are defined in [11,
Theorem 3.1 and Theorem 4.1].

We discuss system (1.2) in $\Omega$. In order to obtain almost
periodic solution and asymptotic stability of system (1.2) we
introduce the following adjoint system (2.1) of system (1.2)
$$
\left\{
    \begin{array}{ll}
\dot{x}_{1}={x_{1}[r_{1}(t)-a_{11}(t)x_{1}-a_{12}(t)y_{1}]+D_{1}(t)(x_{2}-x_{1}),}\vspace{1.5mm} \\
\dot{y}_{1}={y_{1}[r_{2}(t)-a_{21}(t)x_{1}-a_{22}(t)y_{1}]+D_{2}(t)(y_{2}-y_{1}),} \vspace{1.5mm} \\
\dot{x}_{2}={x_{2}[s_{1}(t)-b_{11}(t)x_{2}-b_{12}(t)y_{2}]+D_{1}(t)(x_{1}-x_{2}),} \vspace{1.5mm} \\
\dot{y}_{2}={y_{2}[s_{2}(t)-b_{21}(t)x_{2}-b_{22}(t)y_{2}]+D_{2}(t)(y_{1}-y_{2}),}\vspace{1.5mm} \\
\dot{\tilde{x}}_{1}={\tilde{x}_{1}[r_{1}(t)-a_{11}(t)\tilde{x}_{1}
-a_{12}(t)\tilde{y}_{1}]+D_{1}(t)(\tilde{x}_{2}-\tilde{x}_{1}),}\vspace{1.5mm}  \\
\dot{\tilde{y}}_{1}={\tilde{y}_{1}[r_{2}(t)-a_{21}(t)\tilde{x}_{1}
-a_{22}(t)\tilde{y}_{1}]+D_{2}(t)(\tilde{y}_{2}-\tilde{y}_{1}),} \vspace{1.5mm} \\
\dot{\tilde{x}}_{2}={\tilde{x}_{2}[s_{1}(t)-b_{11}(t)\tilde{x}_{2}
-b_{12}(t)\tilde{y}_{2}]+D_{1}(t)(\tilde{x}_{1}-\tilde{x}_{2}),}\vspace{1.5mm}  \\
\dot{\tilde{y}}_{2}={\tilde{y}_{2}[s_{2}(t)-b_{21}(t)\tilde{x}_{2}
-b_{22}(t)\tilde{y}_{2}]+D_{2}(t)(\tilde{y}_{1}-\tilde{y}_{2}).}
\end{array}
  \right.\eqno(2.1)
$$
Such adjoint system  can be found in [7,12,13].

\section*{ 3. Almost periodic solution and asymptotic stability}
In this section, we will derive some sufficient conditions for
the existence of  almost periodic solution and
its asymptotic stability of system (1.2).\\
\textbf{ Theorem 3.1.}\quad If the conditions of ($H_1$), Lemma 2.2
and Lemma 2.3 hold, and further assume  that system (1.2) satisfies
$$
a^{L}_{11}+a^{L}_{21}>\frac{D^{M}_{1}}{x^L_{2}},\quad
b^{L}_{11}+b^{L}_{21}>\frac{D^{M}_{1}}{x^L_{1}},\quad
a^{L}_{12}+a^{L}_{22}>\frac{D^{M}_{2}}{y^L_{2}},\quad
b^{L}_{12}+b^{L}_{22}>\frac{D^{M}_{2}}{y^L_{1}}.
$$
then system (1.2) has a unique almost periodic
solution which is uniformly asymptotically stable.\\
\textbf{Proof.} Let $z(t)=(x_{1}(t), x_{2}(t), y_{1}(t),y_{2}(t))$
be any solution of system (1.2), from the attraction of $\Omega$, we
discuss our problem in $\Omega$, and denote
$$
Z(t)=(X_{1}(t), X_{2}(t), Y_{1}(t),Y_{2}(t))
$$
We consider the adjoint system (2.1) of system (1.2), let
$${X}_{i}(t)=\ln x_{i}(t),\quad{Y}_{i}(t)=\ln y_{i}(t),
\quad{\tilde{X}}_{i}(t)=\ln
\tilde{x}_{i}(t),\quad{\tilde{Y}}_{i}(t)=\ln
\tilde{y}_{i}(t)(i=1,2),$$ where $x_{1}(t),y_{i}(t),
\tilde{x}_{i}(t), ,\tilde{y}_{i}(t)(i=1,2) $ are the solutions of
adjoint system (2.1) on $\Omega\times\Omega$.

Define Liapunov function:
$$
V(t)=V(t,Z(t),\tilde{Z}(t))=\sum^{2}_{i=1}\mid{X}_{i}(t)-{\tilde{X}}_{i}(t) \mid
+\sum^{2}_{i=1}\mid{Y}_{i}(t)-{\tilde{Y}}_{i}(t) \mid
$$
Taking
$$
a(r)=b(r)=\sum^{2}_{i=1}\mid{X}_{i}(t)-{\tilde{X}}_{i}(t) \mid
+\sum^{2}_{i=1}\mid{Y}_{i}(t)-{\tilde{Y}}_{i}(t) \mid
$$
and $a(r), b(r)$ are continuous increasing positive functions, then
$V(t)$ satisfies the condition (i) of Lemma 2.3. Again from
$$
\begin{array}{ll}
\displaystyle\sum^{2}_{i=1}\mid{X}_{i1}(t)-{\tilde{X}}_{i1}(t) \mid
+\sum^{2}_{i=1}\mid{Y}_{i1}(t)-{\tilde{Y}}_{i1}(t) \mid
-\bigg(\sum^{2}_{i=1}\mid{X}_{i2}(t)-{\tilde{X}}_{i2}(t) \mid
+\sum^{2}_{i=1}\mid{Y}_{i2}(t)-{\tilde{Y}}_{i2}(t)\mid\bigg)\\
\leq \displaystyle\sum^{2}_{i=1}\mid{X}_{i1}(t)-{X}_{i2}(t) \mid
+\sum^{2}_{i=1}\mid{Y}_{i1}(t)-{Y}_{i2}(t) \mid
+\sum^{2}_{i=1}\mid{\tilde{X}}_{i1}(t)-{\tilde{X}}_{i2}(t) \mid
+\sum^{2}_{i=1}\mid{\tilde{Y}}_{i1}(t)-{\tilde{Y}}_{i2}(t)\mid
\end{array}
$$
hence $V(t)$ satisfies the condition (ii) of Lemma 2.3. To check the
condition (iii) of Lemma 2.3, we need to calculate upper-right
derivative of system (2.1), for convenience of statements we denote
$$
A_{i}=sign\big(\mid{X}_{i}(t)-{\tilde{X}}_{i}(t) \mid \big) ,
 B_{i}=sign\big(\mid{Y}_{i}(t)-{\tilde{Y}}_{i}(t) \mid \big),(i=1,2)
$$
\begin{eqnarray}
D^{+}V(t)&=&D^{+}\bigg(\sum^{2}_{i=1}\mid{X}_{i}(t)-{\tilde{X}}_{i}(t) \mid
+\sum^{2}_{i=1}\mid{Y}_{i}(t)-{\tilde{Y}}_{i}(t) \mid\bigg)\nonumber \\
&=&\sum^{2}_{i=1}D^{+}\big(\mid{X}_{i}(t)-{\tilde{X}}_{i}(t) \mid\big)
+\sum^{2}_{i=1}D^{+}\big(\mid{Y}_{i}(t)-{\tilde{Y}}_{i}(t) \mid\big)\nonumber \\
&=&\sum^{2}_{i=1}A_{i}\big({\dot{X}}_{i}(t)-{\dot{\tilde{X}}}_{i}(t) \big)
+\sum^{2}_{i=1}B_{i}\big({\dot{Y}}_{i}(t)-{\dot{\tilde{Y}}}_{i}(t) \big)\nonumber \\
&=&\sum^{2}_{i=1}A_{i}\bigg(\frac{\dot{x}_{i}(t)}{x_{i}(t)}-
\frac{\dot{\tilde{x}}_{i}(t)}{\tilde{x}_{i}(t)} \bigg)+
\sum^{2}_{i=1}B_{i}\bigg(\frac{\dot{y}_{i}(t)}{y_{i}(t)}
-\frac{\dot{\tilde{y}}_{i}(t)}{\tilde{y}_{i}(t)} \bigg)\nonumber \\
&=&-A_{1}a_{11}(t)(x_{1}-\tilde{x}_{1})-A_{1}a_{12}(t)(y_{1}-\tilde{y}_{1})+
A_{1}D_{1}(t)\bigg(\frac{x_{2}}{x_{1}}
-\frac{\tilde{x}_{2} }{\tilde{x}_{1} } \bigg)\nonumber \\
&&-A_{2}b_{11}(t)(x_{2}-\tilde{x}_{2})-A_{2}b_{12}(t)(y_{2}-\tilde{y}_{2})+
A_{2}D_{1}(t)\bigg(\frac{x_{1}}{x_{2}}
-\frac{\tilde{x}_{1} }{\tilde{x}_{2} } \bigg)\nonumber \\
&&-B_{1}a_{21}(t)(x_{1}-\tilde{x}_{1})-B_{1}a_{22}(t)(y_{1}-\tilde{y}_{1})+
B_{1}D_{2}(t)\bigg(\frac{y_{2}}{y_{1}}
-\frac{\tilde{y}_{2} }{\tilde{y}_{1} } \bigg)\nonumber \\
&&-B_{2}b_{21}(t)(x_{2}-\tilde{x}_{2})-B_{2}b_{22}(t)(y_{2}-\tilde{y}_{2})+
B_{2}D_{2}(t)\bigg(\frac{y_{1}}{y_{2}}
-\frac{\tilde{y}_{1} }{\tilde{y}_{2} } \bigg)\nonumber \\
&\leq&-(a^{L}_{11}+a^{L}_{21})\mid x_{1}- \tilde{x}_{1} \mid - (b^{L}_{11}+b^{L}_{21})\mid  x_{2}-
 \tilde{x}_{2}\mid\nonumber \\
&-&(a^{L}_{12}+a^{L}_{22})\mid  y_{1}- \tilde{y}_{1} \mid
-(b^{L}_{12}+b^{L}_{22})\mid  y_{2}- \tilde{y}_{2} \mid +\sum^{4}_{i=1}\tilde{D}_{i}(t) \nonumber
\end{eqnarray}
where
$$\tilde{D}_{1}(t)=\left\{
 \begin{array}{ll}
D_{1}(t)\bigg(\displaystyle\frac{x_{1}}{x_{2}}-\frac{\tilde{x}_{1} }{\tilde{x}_{2} } \bigg)
\qquad X_{2}(t)-\tilde{X}_{2}(t)\geq0 \\
D_{1}(t)\bigg(\displaystyle\frac{\tilde{x}_{1}}{\tilde{x}_{2}}-\frac{x_{1} }{x_{2} } \bigg)
\qquad X_{2}(t)-\tilde{X}_{2}(t)<0
\end{array}
\right.
$$
$$\tilde{D}_{2}(t)=\left\{
 \begin{array}{ll}
D_{1}(t)\bigg(\displaystyle\frac{x_{2}}{x_{1}}-\frac{\tilde{x}_{2} }{\tilde{x}_{1} } \bigg)
\qquad X_{1}(t)-\tilde{X}_{1}(t)\geq0 \\
D_{1}(t)\bigg(\displaystyle\frac{\tilde{x}_{2}}{\tilde{x}_{1}}-\frac{x_{2} }{x_{1} } \bigg)
\qquad X_{1}(t)-\tilde{X}_{1}(t)<0
\end{array}
\right.
$$
$$\tilde{D}_{3}(t)=\left\{
 \begin{array}{ll}
D_{2}(t)\bigg(\displaystyle\frac{y_{1}}{y_{2}}-\frac{\tilde{y}_{1} }{\tilde{y}_{2} } \bigg)
\qquad Y_{2}(t)-\tilde{Y}_{2}(t)\geq0 \\
D_{2}(t)\bigg(\displaystyle\frac{\tilde{y}_{1}}{\tilde{y}_{2}}-\frac{y_{1} }{y_{2} } \bigg)
\qquad Y_{2}(t)-\tilde{Y}_{2}(t)<0
\end{array}
\right.
$$
$$\tilde{D}_{4}(t)=\left\{
 \begin{array}{ll}
D_{2}(t)\bigg(\displaystyle\frac{y_{2}}{y_{1}}-\frac{\tilde{y}_{2} }{\tilde{y}_{1} } \bigg)
\qquad Y_{1}(t)-\tilde{Y}_{1}(t)\geq0 \\
D_{2}(t)\bigg(\displaystyle\frac{\tilde{y}_{2}}{\tilde{y}_{1}}-\frac{y_{2} }{y_{1} } \bigg)
\qquad Y_{1}(t)-\tilde{Y}_{1}(t)<0
\end{array}
\right.
$$
There are the following three cases to consider for $\tilde{D}_{1}(t)$: \\
(i) If $X_{2}(t)>\tilde{X}_{2}(t)$ and $t\geq t^{*}$, then
$$
\tilde{D}_{1}(t)\leq\frac{D_{1}(t)}{x_{2}(t)}(x_{1}(t)
-\tilde{x}_{1}(t))\leq\frac{D^{M}_{1}}{x^L_{2}}|x_{1}(t)-\tilde{x}_{1}(t)|
 $$
(ii) If $X_{2}(t)<\tilde{X}_{2}(t)$ and $t\geq t^{*}$ , then
$$
\tilde{D}_{1}(t)\leq\frac{D_{1}(t)}{\tilde{x}_{2}(t)}(\tilde{x}_{1}(t)
-x_{1}(t))\leq\frac{D^{M}_{1}}{x^L_{2}}|x_{1}(t)-\tilde{x}_{1}(t)|
 $$

(iii) If $X_{2}(t)=\tilde{X}_{2}(t)$, similar to the above analysis, we can get the same result
 as (i)and(ii).
From (i)-(iii), we have
$$
\tilde{D}_{1}(t)\leq\frac{D^{M}_{1}}{x^L_{2}}|x_{1}(t)-\tilde{x}_{1}(t)| \quad for \quad t\geq t^{*}
$$
Considering $\tilde{D}_{2}(t),\tilde{D}_{3}(t),\tilde{D}_{4}(t) $ in the same way we can obtain

$$
\tilde{D}_{2}(t)\leq\frac{D^{M}_{1}}{x^L_{1}}|x_{2}(t)-\tilde{x}_{2}(t)|,\
\tilde{D}_{3}(t)\leq\frac{D^{M}_{2}}{y^L_{2}}|y_{1}(t)-\tilde{y}_{1}(t)|,\
 \tilde{D}_{4}(t)\leq\frac{D^{M}_{2}}{y^L_{1}}|y_{2}(t)-\tilde{y}_{2}(t)|,
 $$
 for $t\geq t^{*}$. It then yields that

$$
\begin{array}{cl}
D^{+}V(t)\leq&\displaystyle-(a^{L}_{11}+a^{L}_{21}-\frac{D^{M}_{1}}{x^L_{2}})\mid x_{1}
- \tilde{x}_{1} \mid - (b^{L}_{11}+b^{L}_{21}-\frac{D^{M}_{1}}{x^L_{1}})\mid  x_{2}-
 \tilde{x}_{2}\mid\vspace {2mm}\\
&\displaystyle-(a^{L}_{12}+a^{L}_{22}-\frac{D^{M}_{2}}{y^L_{2}})\mid  y_{1}- \tilde{y}_{1} \mid
-(b^{L}_{12}+b^{L}_{22}-\frac{D^{M}_{2}}{y^L_{1}})\mid  y_{2}- \tilde{y}_{2} \mid
\end{array}
$$
According to the condition of Theorem 3.1, now we let
$$
P_{1}:=a^{L}_{11}+a^{L}_{21}-\frac{D^{M}_{1}}{x^L_{2}},\quad P_{2}:=b^{L}_{11}
+b^{L}_{21}-\frac{D^{M}_{1}}{x^L_{1}},
$$
$$
P_{3}:=a^{L}_{12}+a^{L}_{22}-\frac{D^{M}_{2}}{y^L_{2}},\quad P_{4}:=b^{L}_{12}
+b^{L}_{22}-\frac{D^{M}_{2}}{y^L_{1}},
$$
and $ \eta=\min\{P_{1},P_{2},P_{3},P_{4}\}>0$
then we get that
$$
D^{+}V(t)\leq-\eta\big(\sum^2_{i=1}|x_{i}-\tilde{x}_{i}|+\sum^2_{i=1}
|y_{i}-\tilde{y}_{i}|\big) \eqno(3.1)
$$

By Mean Value Theorem,we have following formula:
$$
\begin{array}{cl}
|x_{i}- \tilde{x}_{i}|=&\displaystyle|e^{X_{i}}- e^{\tilde{X}_{i}}|=
\zeta_i(t)|X_{i}-\tilde{X}_{i}|\geq x_{i}^L|X_{i}- \tilde{X}_{i}|,\vspace{2mm}\\
|y_{i}- \tilde{y}_{i}|=&\displaystyle|e^{Y_{i}}- e^{\tilde{Y}_{i}}|=
\xi_i(t)|Y_{i}-\tilde{Y}_{i}|\geq y_{i}^L|Y_{i}- \tilde{Y}_{i}|,\ (i=1,2)
\end{array}
$$
where $\zeta_i(t)\in(x_i(t), \tilde{x}_i(t))(i = 1,2)$
and $\xi_i(t)\in(y_i(t), \tilde{y}_i(t))(i = 1,2)$
respectively, then  $\zeta_i(t)\in\Omega$ and
$\xi_i(t)\in\Omega$.

By the above formula and we take $
c=\min\{x_{1}^L\eta,y_{1}^L\eta,x_{2}^L\eta,y_{2}^L\eta\}>0$,  then
we get that
$$
D^{+}V(t)\leq -c\bigg(\sum^{2}_{i=1}\mid{X}_{i}(t)-{\tilde{X}}_{i}(t)
\mid+\sum^{2}_{i=1}\mid{Y}_{i}(t)-{\tilde{Y}}_{i}(t) \mid\bigg)=-c V(t)
$$
It means that $V(t)$ satisfies the condition (iii) of Lemma 2.3. By
Lemma 2.3, system (1.2) has a unique almost periodic solution $z(t)$
on the region $\Omega$, which is uniformly asymptotically stable on
compact set $\Omega$. Since $\Omega$ is the ultimately bounded
region and compact set of system (1.2), hence we get that the
solution $z(t)$ is ultimately bounded on $\Omega$, therefore when
the conditions of Lemma 2.3 hold, almost periodic solution $z(t)$ is
uniformly asymptotically stable. It shows that system (1.2) has a
unique almost periodic solution,
which is uniformly asymptotically stable. This completes the proof.\\

\subsection*{4. Global attractivity}
In this section, we will derive some sufficient conditions for the
globally attractivity of system (1.2).

\textbf{Theorem 4.1}. If the system (1.2) satisfies all the
conditions of Theorem 3.1, then the unique
almost periodic solution of the system (1.2) is globally attractive.\\
\textbf{Proof.} Let $ z(t)=(x_{1}(t), x_{2}(t), y_{1}(t),y_{2}(t))$
be a definitive almost periodic solution of the system (1.2),
$\tilde{z}(t)=(\tilde{x}_{1}(t), \tilde{x}_{2}(t),
\tilde{y}_{1}(t),\tilde{y}_{2}(t))$ be any solution of the system
(1.2).

Construct the same Liapunov function as defined in the proof of
Theorem 3.1,
$$
V(t)=V(t,Z(t),\tilde{Z}(t))=\sum^{2}_{i=1}\mid{X}_{i}(t)-{\tilde{X}}_{i}(t)
\mid +\sum^{2}_{i=1}\mid{Y}_{i}(t)-{\tilde{Y}}_{i}(t) \mid\eqno(4.1)
$$
Integrating both sides of (3.1) from 0 to t, we can derive
$$
V(t)+\eta\int^{t}_{0}\bigg(\sum^{2}_{i=1}\mid{x}_{i}(s)-{\tilde{x}}_{i}(s)
\mid+\sum^{2}_{i=1} \mid{y}_{i}(s)-{\tilde{y}}_{i}(s)
\mid\bigg)ds\leq V(0)\eqno(4.2)
 $$
The expression (4.2) shows that
$$
0\leq V(t)\leq V(0)=\sum^{2}_{i=1}\mid{X}_{i}(0)-{\tilde{X}}_{i}(0)
\mid +\sum^{2}_{i=1}\mid{Y}_{i}(0)-{\tilde{Y}}_{i}(0)
\mid<+\infty,\quad t\geq0  \eqno(4.3)
$$
and
$$
\int^{t}_{0}\bigg(\sum^{2}_{i=1}\mid{x}_{i}(s)-{\tilde{x}}_{i}(s)
\mid+\sum^{2}_{i=1} \mid{y}_{i}(s)-{\tilde{y}}_{i}(s)
\mid\bigg)ds\leq\frac{V(0)}{\eta}<+\infty,\quad t\geq0\eqno(4.4)
 $$
The expression (4.4) implies that

$$
\sum^{2}_{i=1}\mid{x}_{i}(s)-{\tilde{x}}_{i}(s) \mid+\sum^{2}_{i=1}\mid{y}_{i}(s)-{\tilde{y}}_{i}(s)
 \mid\in L^{1}[0,+\infty)\eqno(4.5)
 $$
Obviously $x_i(t)$ and $y_i(t) (i = 1, 2)$ are uniformly bounded, so $X_i(t)$ and $Y_i(t) (i =1, 2)$
 are also uniformly bounded. In addition, by (4.1)-(4.3), we can know that $\tilde{X}_i(t)$ and $\tilde{Y}_i(t) (i =1, 2)$
are uniformly bounded, so $\tilde{x}_i(t)$ and $\tilde{y}_i(t) (i =
1, 2)$ are also uniformly bounded. Combining this fact with the
system (1.2), we have вл$
\dot{x}_i,\dot{y}_i,\dot{\tilde{x}}_i,\dot{\tilde{y}}_i(i = 1, 2)$
are uniformly bounded. Therefore we can easily check
$[x_{i}(t)-\tilde{x}_{i}(t)]$ and $[y_{i}(t)-\tilde{y}_{i}(t)](i=1,
2)$ and their derivatives remain bounded on $[0,+\infty)$. As a
consequence $\sum^{2}_{i=1}\mid{x}_{i}(t)-{\tilde{x}}_{i}(t)
 \mid+\sum^{2}_{i=1}\mid{y}_{i}(t)-{\tilde{y}}_{i}(t) \mid$ is
uniformly continuous on $[0,+\infty)$. From the expression (4.5), it
follows that $\sum^{2}_{i=1}\mid{x}_{i}(t)-{\tilde{x}}_{i}(t)
 \mid+\sum^{2}_{i=1}\mid{y}_{i}(t)-{\tilde{y}}_{i}(t) \mid$ is integrable on $[0,+\infty)$.
By Lemma 2.4, it follows that
$$\lim_{t\rightarrow\infty}\bigg(\sum^{2}_{i=1}\mid{x}_{i}(t)
-{\tilde{x}}_{i}(t) \mid+\sum^{2}_{i=1}\mid{y}_{i}(t)-{\tilde{y}}_{i}(t) \mid\bigg)=0$$
Hence
$$\lim_{t\rightarrow\infty}\mid{x}_{i}(t)-{\tilde{x}}_{i}(t)
\mid=0,\quad\lim_{t\rightarrow\infty}\mid{y}_{i}(t)-{\tilde{y}}_{i}(t)
\mid=0 \quad (i=1,2) $$ This result implies that the unique almost
periodic solution of system (1.2) is stable and attracts all
positive solution of system (1.2). This completes the proof.

\subsection*{5. One example}
Example. we consider the following system
$$
\begin{array}{rl}
\dot{x}_{1}(t)=&\displaystyle
x_{1}(t)\big(5+0.5(sin(\sqrt{2}t)+sin(t))-(2.5+0.5(cos(\sqrt{2}t)+cos(t)))x_{1}(t)
\vspace{1mm}\\
&\displaystyle-(2.2+0.3(sin(\sqrt{2}t)+sin(t)))y_{1}(t)\big)\vspace{1mm} \\
&+(1+0.1(cos(\sqrt{2}t)+cos(t)))(x_{2}(t)-x_{1}(t)),\vspace{1mm} \\
\dot{y}_{1}(t)=&\displaystyle
y_{1}(t)\big(5+0.4(sin(\sqrt{2}t)+sin(t))-(2.25+0.6(cos(\sqrt{2}t)+cos(t)))x_{1}(t)
\vspace{1mm}\\
&\displaystyle-(2.4+0.4(sin(\sqrt{2}t)+sin(t)))y_{1}(t)\big)\vspace{1mm} \\
&+(1+0.2(sin(\sqrt{2}t)+sin(t)))(y_{2}(t)-y_{1}(t)),\vspace{1mm} \\
\dot{x}_{2}(t)=&\displaystyle
x_{2}(t)\big(4+0.5(cos(\sqrt{2}t)+cos(t))-(2.4+0.7(sin(\sqrt{2}t)+sin(t)))x_{2}(t)
\vspace{1mm}\\
&\displaystyle-(2.3+0.5(cos(\sqrt{2}t)+cos(t)))y_{2}(t)\big)\vspace{1mm} \\
&+(1+0.1(cos(\sqrt{2}t)+cos(t)))(x_{2}(t)-x_{1}(t)),\vspace{1mm} \\
\dot{y}_{2}(t)=&\displaystyle
y_{2}(t)\big(4+0.3(cos(\sqrt{2}t)+cos(t))-(2.3+0.5(sin(\sqrt{2}t)+sin(t)))x_{2}(t)
\vspace{1mm}\\
&\displaystyle-(2.5+0.3(cos(\sqrt{2}t)+cos(t)))y_{2}(t)\big)\vspace{1mm} \\
&+(1+0.2(sin(\sqrt{2}t)+sin(t)))(y_{1}(t)-y_{2}(t)),\vspace{1mm} \\
\end{array}
\eqno(5.1)
$$
It is easy to verify all the conditions required in Theorem 3.1 and
Theorem 4.1 are satisfied. Then the system (5.1) has a unique almost
periodic solution and which is  asymptotically stable and globally
attractive.
\begin{center}
    \includegraphics[width=.900\textwidth,height=83mm]{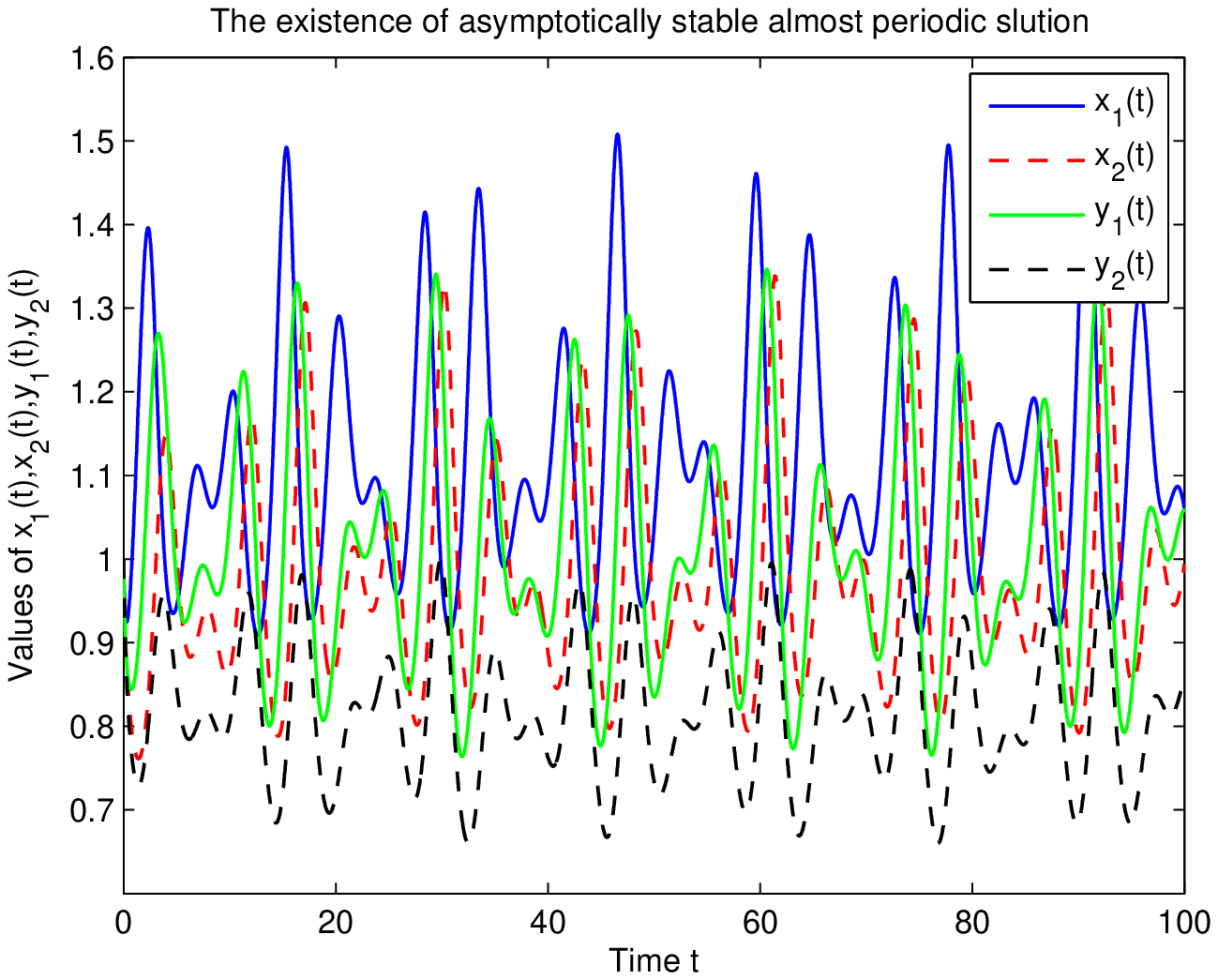}
\put(-435,-10){{\scriptsize{Fig. 1. The existence of almost periodic
solutions for system (5.1). Here, we take the initial value $x_0 =
(x_{10}, x_{20}, y_{10},y_{20}) = (1, 1, 1,1).$}}}
\end{center}
\begin{center}
    \includegraphics[width=.533\textwidth,height=63mm]{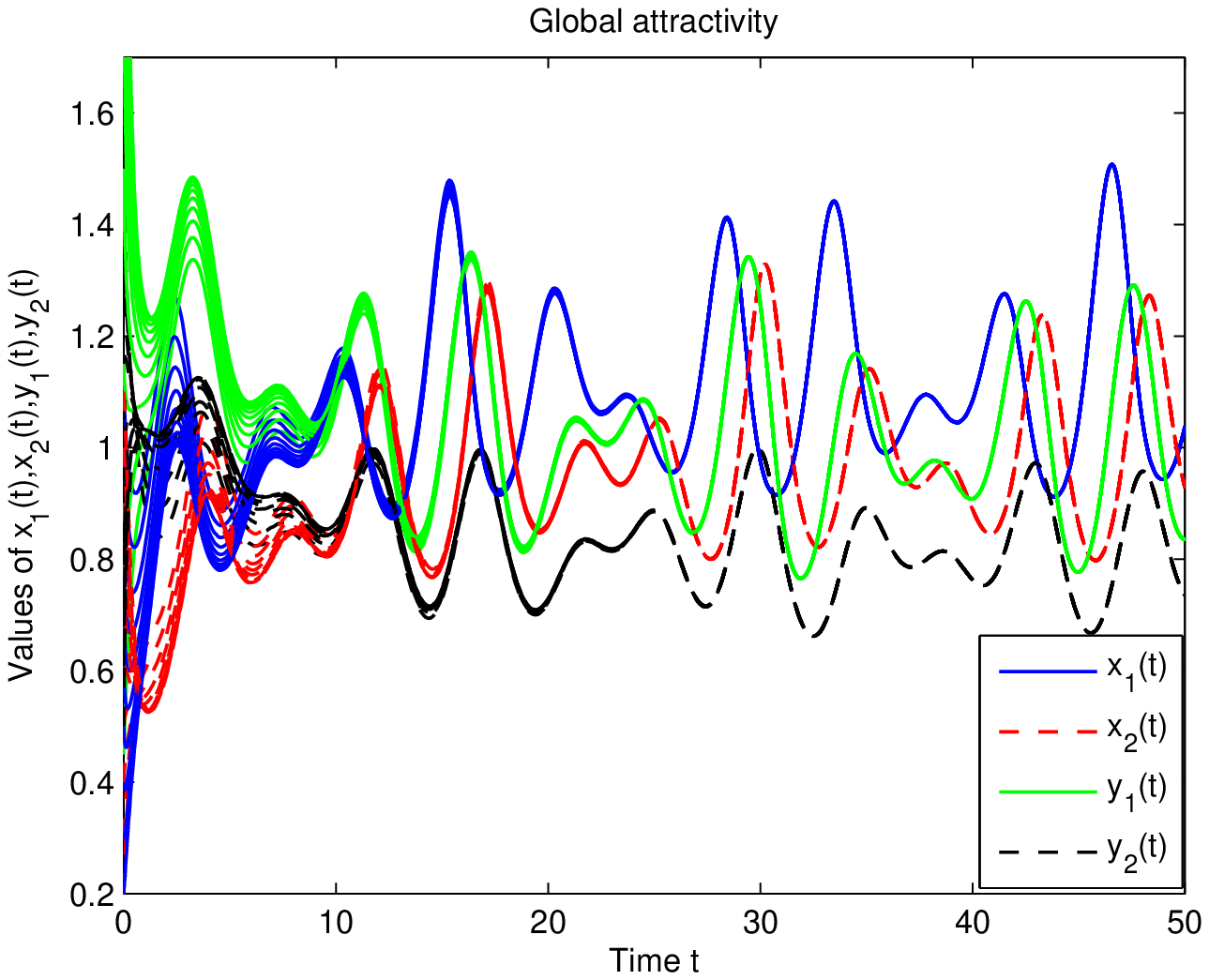}
\put(-205,-10){}
    \includegraphics[width=.53\textwidth,height=63mm]{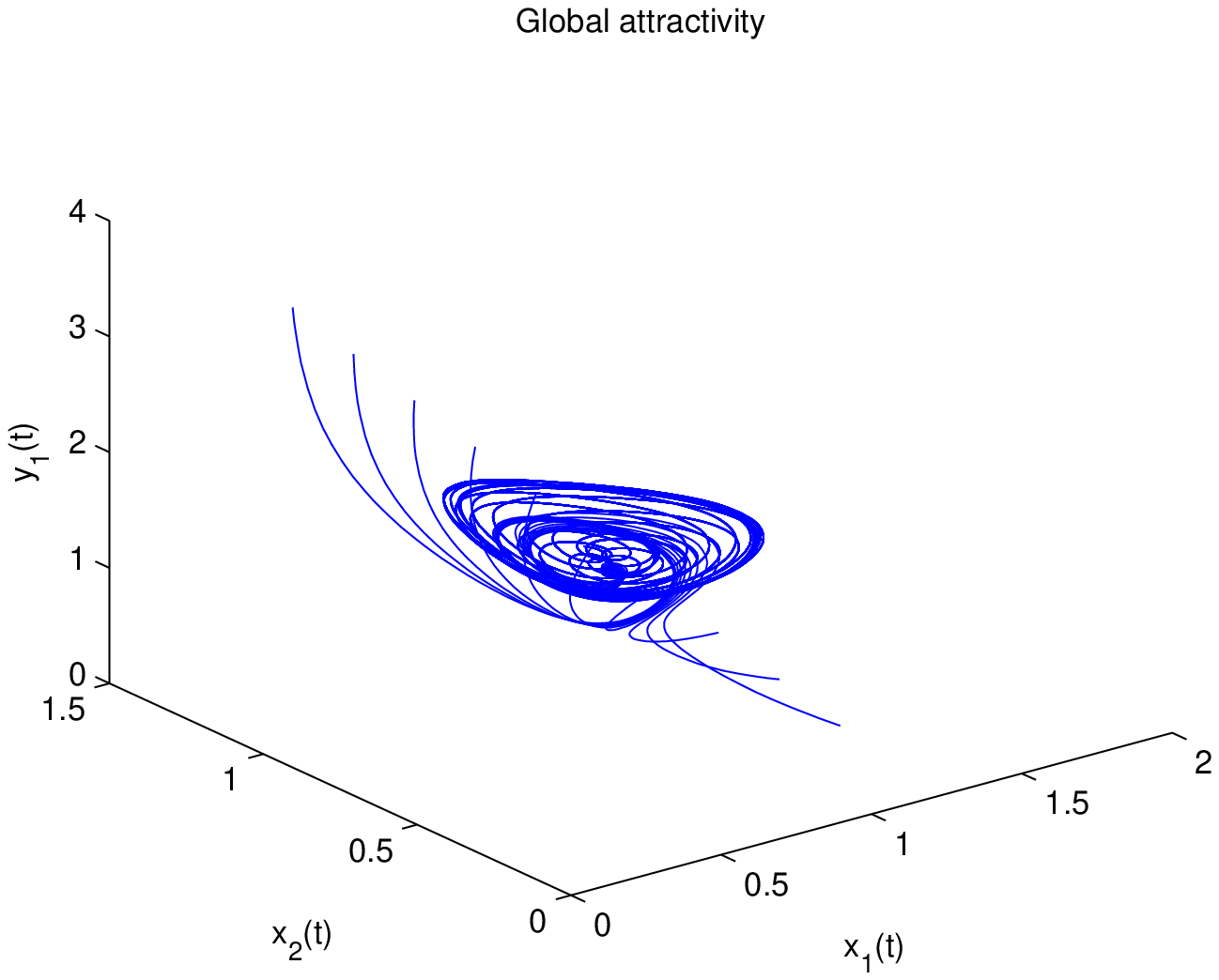}
\put(-205,-10){}
\end{center}

\begin{center}
    \includegraphics[width=.533\textwidth,height=63mm]{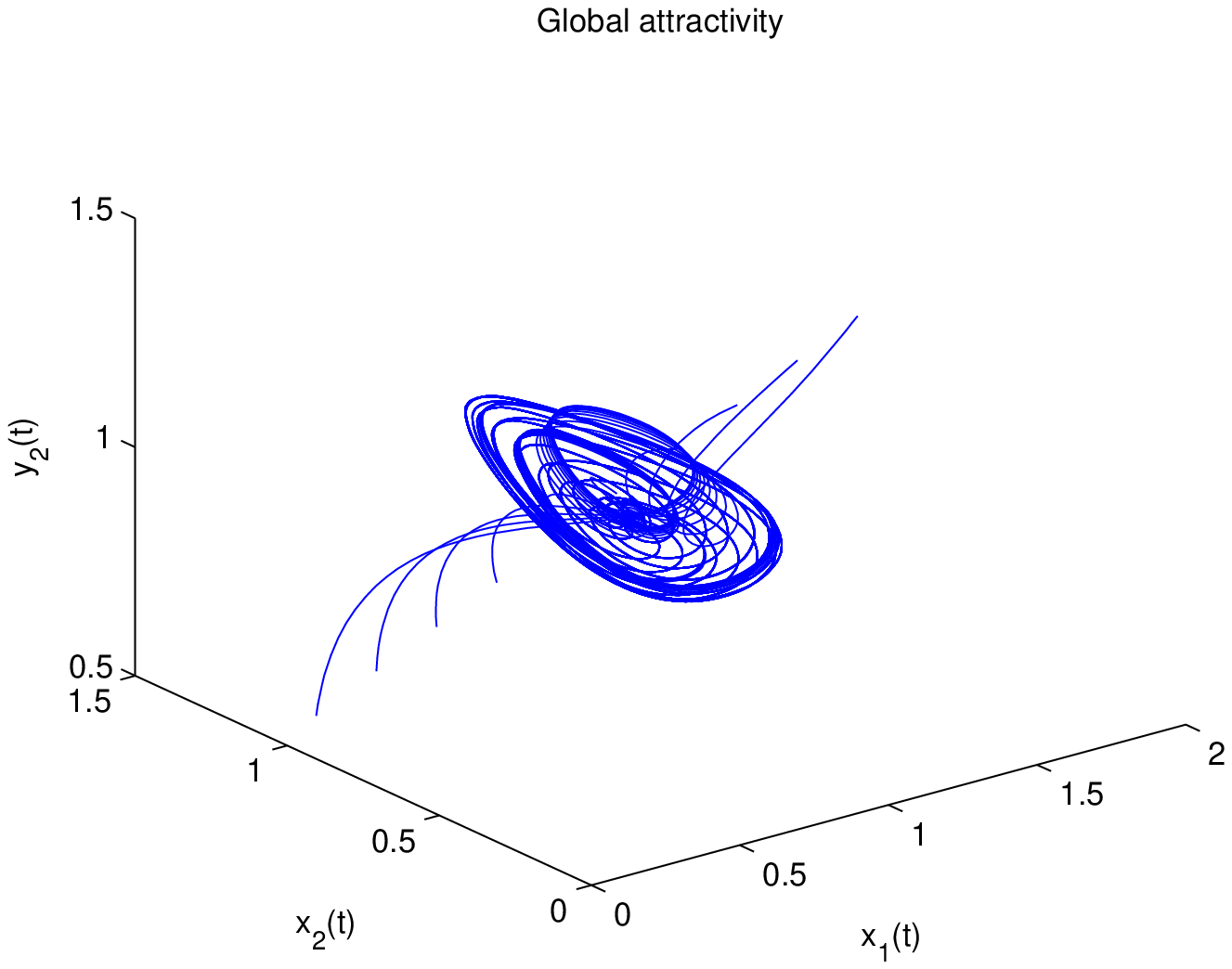}
\put(-205,-10){{\scriptsize{Fig. 2. Global attractivity of almost
periodic solutions for system (5.1). Here, we take different initial
values.}}}
     \includegraphics[width=.53\textwidth,height=63mm]{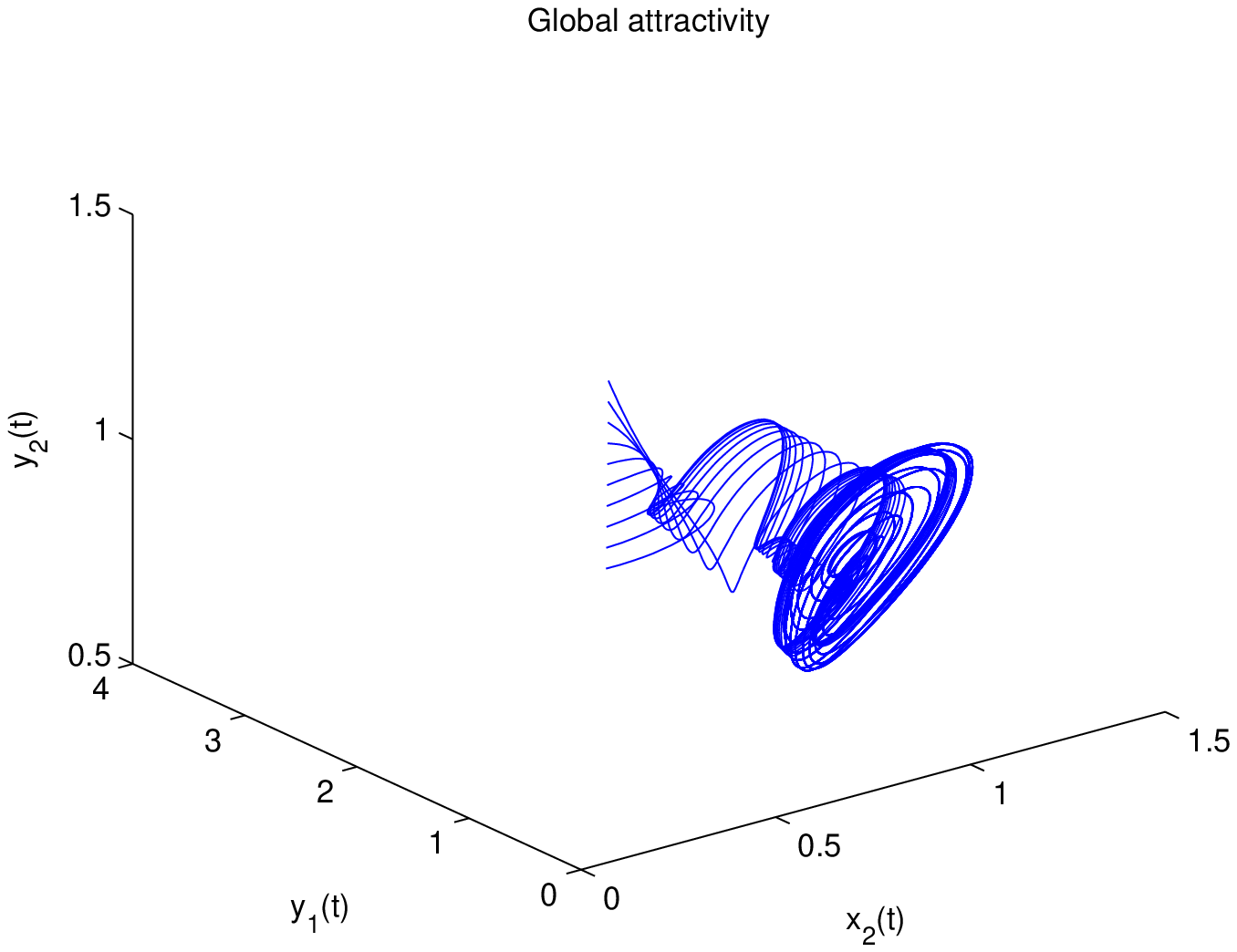}
\put(-215,-10){}
\end{center}

Figure 1 shows that the system (5.1) converges to a
 almost periodic solution. Figure 2  shows that the almost periodic solution
  of system (5.1) is globally attractive.

\def\refname

\end{document}